\documentclass[11pt]{amsart}
\usepackage{amssymb}
\usepackage{amsmath}
\usepackage{mathtools, float}
\usepackage{amsfonts}
\usepackage{stmaryrd}
\usepackage[english]{babel}
\usepackage[utf8]{inputenc}
\usepackage{enumerate}
\usepackage{graphicx,pinlabel,xcolor, fullpage}
\usepackage{adjustbox}
\usepackage{array,multirow}
\usepackage{float}

\theoremstyle{plain}

\newtheorem{theorem*}{Theorem}

\theoremstyle{definition}

\newtheoremstyle{case}{}{}{}{}{}{:}{ }{}
\theoremstyle{case}

\newcommand{\QED}{\hfill $\blacksquare$}

\begin{document}
\title[On the Schur Multiplier of finite $p$-groups of maximal class]{On the Schur Multiplier of finite $p$-groups of maximal class}
\author{Renu Joshi and Siddhartha Sarkar}
\address{Department of Mathematics\\
Indian Institute of Science Education and Research Bhopal\\
Bhopal Bypass Road, Bhauri \\
Bhopal 462 066, Madhya Pradesh\\
India}
\email{renu16@iiserb.ac.in, sidhu@iiserb.ac.in}
%\urladdr{https://home.iiserb.ac.in/~sidhu/}

%\keywords{Finite $p$-groups; Genus spectrum; $p$-groups of maximal class}

\subjclass{
Primary 20D15, Secondary 20J06
}

\begin{abstract}
In this article, we prove that the Schur Multiplier of a finite $p$-group of maximal class of order $p^n ~(4 \leq n \leq p+1)$ is elementary abelian. The case $n = p+1$ settles a question raised by Primo\v{z} Moravec in an earlier article.        
\end{abstract}

\maketitle

\section{Introduction}
\label{introsec}

\noindent Let $G$ be a finite group. The Schur multiplier $M(G)$ of $G$ is defined to be the second integral homology group $H_2(G, {\mathbb{Z}})$, where ${\mathbb{Z}}$ is considered to be the trivial $G$-module. This was introduced by Schur \cite{sch} to study projective representations of finite groups. In the last century, it was used to study central group extensions and consequently plays a central role in studying the classification of finite $p$-groups. Several authors have determined the upper bound of the order, rank, and exponent of $M(G)$ for finite groups (see \cite{bet2, mor}). As a particularly interesting case, the finite $p$-groups with trivial Schur multipliers were studied in \cite{bet}.

\smallskip

\noindent A constructive approach to defining the Schur multiplier is related to constructing the non-abelian tensor square $G \otimes G$ of $G$. This is defined as follows: the generators of $G \otimes G$ are the abstract symbols $g \otimes h$ for all $g, h \in G$ and the relations of this group are the relations:
\begin{eqnarray*}
(g_1 g_2) \otimes h & = & (g^{g_2}_1 \otimes h^{g_2})(g_2 \otimes h) \\
g \otimes (h_1 h_2) & = & (g \otimes h_2)(g^{h_2} \otimes h^{h_2}_1)
\end{eqnarray*}  
for all $g, g_1, g_2, h, h_1, h_2 \in G$ and $x^y = y^{-1}xy$ for all $x, y \in G$. Let $\nabla(G)$ denote the normal subgroup of $G$ generated by the elements $g \otimes g$ for all $g \in G$. The exterior square $G \wedge G$ is defined to be the quotient group $(G \otimes G)/{\nabla(G)}$. The function $f : G \wedge G \rightarrow [G, G]$ defined on the generators by $f(g \wedge h) := [g,h] = g^{-1} h^{-1} gh ~(g, h \in G)$ is a group epimorphism whose kernel is isomorphic to $M(G)$ (see \cite{bl}).
 
\smallskip

\noindent In this article, we consider the construction of a group $\nu(G)$ due to Rocco \cite{roc}. Let $G$ be a group and $G^{\varphi}$ be an isomorphic copy of $G$, where $\varphi$ is an isomorphism. Let ${\mathcal{N}}$ denote the normal subgroup of the free product $G \ast G^{\varphi}$ generated by the relations:
\[
[g_1, g^{\varphi}_2]^{g_3} = [g^{g_3}_1, (g^{g_3}_2)^{\varphi}] = [g_1, g^{\varphi}_2]^{g^{\varphi}_3} \hspace*{.2in} (g_1, g_2, g_3 \in G).
\] 
Then $\nu(G)$ is defined to be the quotient $(G \ast G^{\varphi})/{\mathcal{N}}$. If $G$ is a finite group (resp. a finite $p$-group of nilpotency class $c$), then $\nu(G)$ is a finite group (resp. a finite $p$-group of nilpotency class at most $c+1$) (see \cite[Proposition 2.4 and Corollary 3.2]{roc}). 

\smallskip

\noindent For subgroups $H, K \leq G$, we define
\[
[H, K^{\varphi}] := \langle [a, b^{\varphi}] ~:~ a \in H, b \in K \rangle. 
\]
Then for any normal subgroup $N \unlhd G$, we have $[N, G^{\varphi}]$ and $[G, N^{\varphi}]$ are normal subgroups of $\nu(G)$ (\cite[Proposition 2.5]{roc}). The subgroup $\kappa(G) \leq [G, G^{\varphi}]$ defined as $\kappa(G) := \langle [g, g^{\varphi}] ~:~ g \in G \rangle$ is a normal subgroup of $[G, G^{\varphi}]$. The function $T : G \otimes G \rightarrow [G, G^{\varphi}]$ defined on the generators of $G \otimes G$ as 
\[
T(a \otimes b) := [a, b^{\varphi}]
\]  
induces an isomorphism of groups. This induces an isomorphism ${\overline{T}} : G \wedge G \rightarrow [G, G^{\varphi}]/{\kappa(G)}$ defined as ${\overline{T}}(a \wedge b) = [a, b^{\varphi}] \kappa(G) =: [[a, b^{\varphi}]]$. So we can identify the groups $G \wedge G$ and $[[G, G^{\varphi}]] := [G, G^{\varphi}]/{\kappa(G)}$ by ${\overline{T}}$ and regard the Schur Multiplier $M(G)$ of $G$ as the kernel of the map $\Psi : [[G, G^{\varphi}]] \rightarrow [G,G]$ defined as $\Psi([[a, b^{\varphi}]]) := [a,b] ~(a, b \in G)$. The advantage of working in $[G, G^{\varphi}]$ (resp. $[[G, G^{\varphi}]]$) instead of $G \otimes G$ (resp. $G \wedge G$) is that we can apply the commutator identities and the Hall's collection process in $\nu(G)$ which is used very effectively in \cite{bmo, hky}.    

\bigskip

\noindent Let $G$ be a finite $p$-group of nilpotency class $c$. The co-class of $G$ is defined as $\log_p|G| - c$. This article is related to computing Schur multipliers of finite $p$-groups of co-class $1$ (also called the finite $p$-groups of maximal class). If $|G| = p^n$, then we have $n \geq 2$. For $n=2$ we have $M({\mathbb{Z}}_{p^2}) = 1$ and $M({\mathbb{Z}}_p \times {\mathbb{Z}}_p) \cong {\mathbb{Z}}_p$, where ${\mathbb{Z}}_m$ denotes the finite cyclic group of order $m$. For $n = 3$, and $p$ odd, we also know $M(M_{p^3}) = 1$ and $M(H_{p^3}) = {\mathbb{Z}}_p \times {\mathbb{Z}}_p$, where $M_{p^3}$ and $H_{p^3}$ are non-abelian groups of order $p^3$ and exponent $p^2$ and $p$ respectively. For $p=2$, the Schur multipliers of $2$-groups of co-class $1$ are all known: $M({\mathbb{D}}_{2^n}) \cong {\mathbb{Z}}_2, M({\mathbb{Q}}_{2^n}) = 1, M({\mathbb{SD}}_{2^n}) = 1$, where ${\mathbb{D}}_{2^n}$ and ${\mathbb{Q}}_{2^n}$ are dihedral and quaternion groups of order $2^n ~(n \geq 3)$ and ${\mathbb{SD}}_{2^n}$ is Semi-dihedral group of order $2^n ~(n \geq 4)$ (see \cite[Theorem 2.11.3]{kar}). For $p$ odd, the Schur multiplier of all finite $p$-groups of order up to $p^5$ can be found in \cite{hky}.  

\bigskip

\noindent For the remaining part of the paper, we assume $p$ is an odd prime and $G$ is a finite $p$-group of co-class $1$ of order $p^n ~(n \geq 4)$. We denote $P_2(G) := [G, G]$ and $P_i(G) := [P_{i-1}(G), G]$ for $3 \leq i \leq n-1$. The subgroups $P_1(G) := C_G(P_2(G)/{P_4(G)})$ and $C_G(P_{n-2}(G))$ are maximal subgroups of $G$, which are called the two step centralizers of $G$. Let $s_0 \in G \setminus (P_1(G) \cup C_G(P_{n-2}(G)))$ and $s_1 \in P_1(G) \setminus P_2(G)$. Define $s_i := [s_{i-1}, s_0]$ for all $i \geq 2$. Setting $P_0(G) := G$, each sections $P_i(G)/{P_{i+1}(G)}$ are cyclic of order $p$ and are generated by $s_i P_{i+1}(G) ~(i \geq 0)$ (see \cite[Chapter 3]{lmc}).    

\bigskip

\noindent The methods in this article are developed to resolve the question: if $G$ is a finite $p$-group of order $p^{p+1}$ and maximal class, then $M(G)$ is elementary abelian. This question was raised by Primo\v{z} Moravec in \cite{mor}, and it was proved that if $G$ is a finite $p$-group of maximal class, then ${\mathrm{exp}}(M(G)) \leq {\mathrm{exp}}(G)$. In particular, if ${\mathrm{exp}}(G) = p$, then necessarily $n \leq p$, and we have ${\mathrm{exp}}(M(G)) \leq p$. In case ${\mathrm{exp}}(G) = p^2$, then $p^4 \leq |G| \leq p^{2p-1}$ and from the same result it follows that ${\mathrm{exp}}(M(G)) \leq p^2$.   

\bigskip

\noindent The main results of this article are as follows: 

\bigskip

\noindent {\bf Theorem 1.} (Theorem \ref{gen-at-least-wt2}, Corollary \ref{cor-gen-at-least-wt2}) Let $p$ be an odd prime and $G$ be a finite $p$-group of maximal class of order $p^n ~(4 \leq n \leq p+1)$ and exponent $p^2$. Then $[P_i(G), G^{\varphi}] = [P_i(G)^{\varphi}, G] \leq \nu(G)$, and ${\mathrm{exp}}([P_i(G), G^{\varphi}]) \leq p$ for all $2 \leq i \leq n-1$. 

\bigskip

\noindent As a consequence of this, we derive:

\bigskip

\noindent {\bf Theorem 2.} (Theorem \ref{main-thm-Schur-mult}) Let $G$ be a finite $p$-group of maximal class and order $p^n ~(4 \leq n \leq p+1)$. Then the Schur multiplier $M(G)$ of $G$ is elementary abelian.

\bigskip

\noindent This result is the best possible. Denoting $G(n,r) := {\mathrm{SmallGroup}}(n,r)$, here is a complete list of groups $G$ of order $3^5$ and $5^7$ with $M(G)$ not elementary abelian (computed using HAP \cite{hap} package of GAP \cite{gap}):

\bigskip

\begin{center}
\begin{tabular}{ |c|c|c|c| } 
 \hline
 $(n,r)$ & $M(G)$ & $(n,r)$ & $M(G)$ \\ 
 \hline \hline
 $(243,26)$ & ${\mathbb{Z}}_3 \times {\mathbb{Z}}_9$ & $(78125,1283)$ & ${\mathbb{Z}}_5 \times {\mathbb{Z}}_5 \times {\mathbb{Z}}_{25}$ \\ 
 \hline
 $(243,28)$ & ${\mathbb{Z}}_9$ & $(78125,1286)$ & ${\mathbb{Z}}_5 \times {\mathbb{Z}}_5 \times {\mathbb{Z}}_{25}$ \\ 
 \hline
  &  & $(78125,1297)$ & ${\mathbb{Z}}_5 \times {\mathbb{Z}}_5 \times {\mathbb{Z}}_{25}$ \\
 \cline{3-4} 
 &  & $(78125,1304)$ & ${\mathbb{Z}}_5 \times {\mathbb{Z}}_5 \times {\mathbb{Z}}_{25}$ \\
 \cline{3-4} 
 &  & $(78125,1370)$ & ${\mathbb{Z}}_5 \times {\mathbb{Z}}_5 \times {\mathbb{Z}}_{25}$ \\
 \hline
\end{tabular}
\end{center}

\bigskip

\noindent We are informed by Professor Michael VaughanLee that there are $17212$ groups of order $7^9$ of maximal class whose Schur multipliers are isomorphic to ${\mathbb{Z}}^2_7$, ${\mathbb{Z}}^3_7$, ${\mathbb{Z}}^4_7$, ${\mathbb{Z}}_7 \times {\mathbb{Z}}_7 \times {\mathbb{Z}}_7 \times {\mathbb{Z}}_{49}$ and ${\mathbb{Z}}_7 \times {\mathbb{Z}}_7 \times {\mathbb{Z}}_{49}$ and these appeared respectively $12596, 4487, 84, 25$ and $20$ times.  

\bigskip
%\newpage

\section{Preliminaries}
\label{prelim}

\noindent In this section, we will discuss some basic results that will be used later. 

\bigskip

\subsection{Lemma}\label{rels-nu-G} (\cite[Lemma 9]{bmo}) For any group $G$, the following relations hold in $\nu(G)$:

\smallskip

\noindent (i) $[g^{\varphi}_1, g_2, g_3] = [g_1, g^{\varphi}_2, g_3] = [g_1, g_2, g^{\varphi}_3] = [g^{\varphi}_1, g^{\varphi}_2, g_3] = [g^{\varphi}_1, g_2, g^{\varphi}_3] = [g_1, g^{\varphi}_2, g^{\varphi}_3]$ for all $g_1, g_2, g_3 \in G$.

\smallskip

\noindent (ii) $[[g_1, g^{\varphi}_2], [h_1, h^{\varphi}_2]] = [[g_1, g_2], [h_1, h^{\varphi}_2]]$ for all $g_1, g_2, h_1, h_2 \in G$. 

\bigskip

\noindent The following lemma is an extension of Lemma \ref{rels-nu-G}(i).

\subsection{Lemma}\label{exponent-1-phi}
Let $G$ be any group and $r \geq 3$ be an integer. For any distinct ${\underline{\epsilon}} = (\epsilon_1, \dotsc, \epsilon_r), {\underline{\delta}} = (\delta_1, \dotsc, \delta_r) \in \{ 1, \varphi \}^r \setminus \{ (1, \dotsc, 1), (\varphi, \dotsc, \varphi) \}$ and $g_1, \dotsc, g_r \in G$ we have the following relation in $\nu(G)$:
\[
[g^{\epsilon_1}_1, \dotsc, g^{\epsilon_r}_r] = [g^{\delta_1}_1, \dotsc, g^{\delta_r}_r].
\]  

\bigskip

\noindent {\bf Proof.} We prove it by induction on $r \geq 3$ and using Lemma \ref{rels-nu-G}. For $r=3$, it follows from Lemma \ref{rels-nu-G}(i). Assume the result for some $r-1 \geq 3$ and consider tuples ${\underline{\epsilon}}, {\underline{\delta}}$ of length $r$ as above.  

\smallskip

\noindent {\bf Case I.} $(\epsilon_1, \dotsc, \epsilon_{r-1}) = (1, \dotsc, 1), (\delta_1, \dotsc, \delta_{r-1}) = (\varphi, \dotsc, \varphi)$.

\smallskip

\noindent This implies $\epsilon_r = \varphi, \delta_r = 1,$ and we have
\begin{eqnarray*}
[g_1, \dotsc, g_{r-1}, g^{\varphi}_r] & = & \Big[ [g_1, \dotsc, g_{r-2}], g_{r-1}, g^{\varphi}_r \Big] \\
	& = & \Big[ [g_1, \dotsc, g_{r-2}]^{\varphi}, g^{\varphi}_{r-1}, g_r \Big] ~=~ [g^{\varphi}_1, \dotsc, g^{\varphi}_{r-2}, g^{\varphi}_{r-1}, g_r].
\end{eqnarray*} 

\smallskip

\noindent The case $(\epsilon_1, \dotsc, \epsilon_{r-1}) = (\varphi, \dotsc, \varphi), (\delta_1, \dotsc, \delta_{r-1}) = (1, \dotsc, 1)$ is similar.

\smallskip

\noindent {\bf Case II.} $(\epsilon_1, \dotsc, \epsilon_{r-1}) = (1, \dotsc, 1)$ and $(\delta_1, \dotsc, \delta_{r-1}) \in \{ 1, \varphi \}^{r-1} \setminus \{ (1, \dotsc, 1), (\varphi, \dotsc, \varphi) \}$.

\smallskip

\noindent In this case we have $\epsilon_r = \varphi$, $(\delta_1, \dotsc, \delta_{r-1})$ is of length $r-1 \geq 3$, and contains at least one $1$ and another $\varphi$. We have two choices: $\delta_r$ is either $1$ or $\varphi$. In case $\delta_r = 1$, we have
\begin{eqnarray*}
[g_1, \dotsc, g_{r-1}, g^{\varphi}_r] & = & \Big[ [g_1, \dotsc, g_{r-2}], g_{r-1}, g^{\varphi}_r \Big] \\
	& = & \Big[ [g_1, \dotsc, g_{r-2}], g^{\varphi}_{r-1}, g_r \Big] ~=~ \Big[ [g_1, \dotsc, g_{r-2}, g^{\varphi}_{r-1}], g_r \Big] \\
	& = & \Big[ [g^{\delta_1}_1, \dotsc, g^{\delta_{r-2}}_{r-2}, g^{\delta_{r-1}}_{r-1}], g_r \Big] ~=~ [g^{\delta_1}_1, \dotsc, g^{\delta_r}_r].
\end{eqnarray*}

\smallskip

\noindent In case $\delta_r = \varphi$, we have
\begin{eqnarray*}
[g_1, \dotsc, g_{r-1}, g^{\varphi}_r] & = & \Big[ [g_1, \dotsc, g_{r-2}], g_{r-1}, g^{\varphi}_r \Big] \\
	& = & \Big[ [g_1, \dotsc, g_{r-2}], g^{\varphi}_{r-1}, g^{\varphi}_r \Big] ~=~ \Big[ [g_1, \dotsc, g_{r-2}, g^{\varphi}_{r-1}], g^{\varphi}_r \Big] \\
	& = & \Big[ [g^{\delta_1}_1, \dotsc, g^{\delta_{r-2}}_{r-2}, g^{\delta_{r-1}}_{r-1}], g^{\varphi}_r \Big] ~=~ [g^{\delta_1}_1, \dotsc, g^{\delta_r}_r].
\end{eqnarray*}

\smallskip

\noindent The case $(\epsilon_1, \dotsc, \epsilon_{r-1}) = (\varphi, \dotsc, \varphi)$ and $(\delta_1, \dotsc, \delta_{r-1}) \in \{ 1, \varphi \}^{r-1} \setminus \{ (1, \dotsc, 1), (\varphi, \dotsc, \varphi) \}$ is similar.

\smallskip

\noindent {\bf Case III.} $(\epsilon_1, \dotsc, \epsilon_{r-1}), (\delta_1, \dotsc, \delta_{r-1}) \in \{ 1, \varphi \}^{r-1} \setminus \{ (1, \dotsc, 1), (\varphi, \dotsc, \varphi) \}$.

\smallskip

\noindent Using induction hypothesis, we have
\[
[g^{\epsilon_1}_1, \dotsc, g^{\epsilon_{r-1}}_{r-1}] = [g^{\delta_1}_1, \dotsc, g^{\delta_{r-1}}_{r-1}]
\] 
and consequently,
\[
\Big[ [g^{\epsilon_1}_1, \dotsc, g^{\epsilon_{r-1}}_{r-1}], g_r \Big] = \Big[ [g^{\delta_1}_1, \dotsc, g^{\delta_{r-1}}_{r-1}], g_r \Big].
\] 
Using the symmetry of both sides, it is enough to show that
\[
\Big[ [g^{\epsilon_1}_1, \dotsc, g^{\epsilon_{r-1}}_{r-1}], g_r \Big] = \Big[ [g^{\epsilon_1}_1, \dotsc, g^{\epsilon_{r-1}}_{r-1}], g^{\varphi}_r \Big].
\]
Now since $r-3 \geq 1$, we have
\begin{eqnarray*}
\Big[ [g^{\epsilon_1}_1, \dotsc, g^{\epsilon_{r-1}}_{r-1}], g_r \Big] & = & \Big[ [g_1, \dotsc, g_{r-2}, g^{\varphi}_{r-1}], g_r \Big] \\
	& = & \Big[ [g_1, \dotsc, g_{r-2}], g^{\varphi}_{r-1}, g_r \Big] \\
	& = & \Big[ [g_1, \dotsc, g_{r-2}], g^{\varphi}_{r-1}, g^{\varphi}_r \Big] \\
	& = & \Big[ [g_1, \dotsc, g_{r-2}, g^{\varphi}_{r-1}], g^{\varphi}_r \Big] ~=~ \Big[ [g^{\epsilon_1}_1, \dotsc, g^{\epsilon_{r-1}}_{r-1}], g^{\varphi}_r \Big].
\end{eqnarray*}
\QED

\bigskip

\noindent We are going to use a few elementary results several times, which we will mention here.

\bigskip

\subsection{Proposition}\label{exponent-Pi-subgroups} (\cite[Proposition 3.3.2]{lmc}) Let $p$ be an odd prime and $G$ be a finite $p$-group of maximal class of order $p^n ~(4 \leq n \leq p+1)$. Then ${\mathrm{exp}}(P_2(G)) = p = {\mathrm{exp}}(G/P_{n-1}(G))$.

\bigskip

\subsection{Proposition}\label{comm-collection} (\cite[Proposition 1.1.32]{lmc}) Let $G$ be any group, $p$ a prime, and $x, y \in G$. For any $r \in {\mathbb{N}}$, we have:
\[
(xy)^{p^r} \equiv x^{p^r} y^{p^r} [y,x]^{p^r \choose 2} [y,x,x]^{p^r \choose 3} \dotsc [y,_{p^r - 1} x] ~{\mathrm{mod}}~ K(x,y), 
\]
\[
[x^{p^r}, y] \equiv [x,y]^{p^r} [x,y,x]^{p^r \choose 2} \dotsc \big[ [x,y],_{p^r - 1} x \big] ~{\mathrm{mod}}~ K(x,[x,y]),
\]

\noindent where $K(a,b)$ denote the normal subgroup generated by the set of all basic commutators in $\{ a, b \} \subseteq \langle x, y \rangle$ of weight at least $p^r$, and weight at least $2$ in $b$, together with the $p^{r-k+1}$-th power of all basic commutators in $\{ a, b \}$ of weight $< p^k$, and weight at least $2$ in $b$ for $1 \leq k \leq r$.

\bigskip

\section{Proof of main theorems}
\label{main-thms}

\bigskip

\noindent In this section, $p$ always denotes an odd prime. We first need the following result describing certain abelian quotients in the group $[G, G^{\varphi}]$.

\bigskip

\subsection{Lemma}\label{section-lemma} Let $G$ be a finite $p$-group of maximal class of order $p^n$. For any integer $0 \leq i \leq n-1$, we have: 

\smallskip

\noindent (i) $[P_i(G), G^{\varphi}]/[P_{i+1}(G), G^{\varphi}]$ is an abelian group generated by the elements
\[
[s_i, s^{\varphi}_j][P_{i+1}(G), G^{\varphi}] \hspace*{.2in} (0 \leq j \leq n-1).
\]

\smallskip

\noindent (ii) $[P_i(G)^{\varphi}, G]/[P_{i+1}(G)^{\varphi}, G]$ is an abelian group generated by the elements
\[
[s^{\varphi}_i, s_j][P_{i+1}(G)^{\varphi}, G] \hspace*{.2in} (0 \leq j \leq n-1).
\]

\bigskip

\noindent {\bf Proof.} (i) Let $[a, b^{\varphi}]$ be a generator of $[P_i(G), G^{\varphi}]$ with $a \in P_i(G), b \in G$. Write $a = s^{\lambda}_i w$ for some $0 \leq \lambda \leq p-1$ and $w \in P_{i+1}(G)$. Then, using commutator identities we have,
\begin{eqnarray*}
[a, b^{\varphi}] & = & [s^{\lambda}_i, b^{\varphi}][s^{\lambda}_i, b^{\varphi}, w][w, b^{\varphi}] \\
	& = & [s^{\lambda}_i, b^{\varphi}][s^{\lambda}_i, b, w^{\varphi}][w, b^{\varphi}] \equiv [s^{\lambda}_i, b^{\varphi}] ~{\mathrm{mod}}~ [P_{i+1}(G), G^{\varphi}].
\end{eqnarray*}

\smallskip

\noindent If $\lambda \geq 2$, then 
\begin{eqnarray*}
[s^{\lambda}_i, b^{\varphi}] & = & [s^{\lambda - 1}_i, b^{\varphi}][s^{\lambda - 1}_i, b^{\varphi}, s_i][s_i, b^{\varphi}] \\
	& = & [s^{\lambda - 1}_i, b^{\varphi}][s^{\lambda - 1}_i, b, s^{\varphi}_i][s_i, b^{\varphi}] \equiv [s^{\lambda - 1}_i, b^{\varphi}][s_i, b^{\varphi}] ~{\mathrm{mod}}~ [P_{i+1}(G), G^{\varphi}]
\end{eqnarray*}
since $[P_{i+1}(G), G^{\varphi}] \unlhd \nu(G)$. Thus by induction on $\lambda$ we get
\[
[a, b^{\varphi}] \equiv [s_i, b^{\varphi}]^{\lambda} ~{\mathrm{mod}}~ [P_{i+1}(G), G^{\varphi}].
\]

\noindent Now we notice that, for any two integers $j_1, j_2 \geq 0$, we have
\[
\Big[ [s_i, s^{\varphi}_{j_1}], [s_i, s^{\varphi}_{j_2}] \Big] = \Big[ [s_i, s_{j_1}], [s_i, s_{j_2}]^{\varphi} \Big] \in [P_{i+1}(G), G^{\varphi}].
\]
Hence the elements $[s_i, s^{\varphi}_j][P_{i+1}(G), G^{\varphi}] ~(0 \leq j \leq n-1)$ centralize each other. 

\smallskip 

\noindent Finally, if $b = b_1 s_j$ for some $j \geq 0$, then we have
\begin{eqnarray*}
[s_i, b^{\varphi}] & = & [s_i, s^{\varphi}_j][s_i, b^{\varphi}_1][s_i, b^{\varphi}_1, s^{\varphi}_j] ~=~ [s_i, s^{\varphi}_j][s_i, b^{\varphi}_1][s_i, b_1, s^{\varphi}_j] \\
	& \equiv & [s_i, s^{\varphi}_j][s_i, b^{\varphi}_1] ~{\mathrm{mod}}~ [P_{i+1}(G), G^{\varphi}].
\end{eqnarray*}

\smallskip

\noindent Since any $b \in G$ can be written as  
\[
b = s^{\lambda_0}_0 s^{\lambda_1}_1 \dotsc s^{\lambda_{n-1}}_{n-1} \hspace*{.2in} (0 \leq \lambda_0, \lambda_1, \dotsc, \lambda_{n-1} \leq p-1),
\]
from above two statements, the statement (i) follows. Proof of (ii) is similar. \QED

\bigskip

\subsection{Theorem}\label{gen-at-least-wt2} Let $G$ be a finite $p$-group of maximal class and order $p^n ~(4 \leq n \leq p+1)$.
For integers $2 \leq i \leq n-1$ and $0 \leq j \leq n-1$ we have:

\smallskip

\noindent (i) All basic commutators in $\{ s_i, [s_i, s^{\varphi}_j] \}$ of weight $\geq p$ and weight at least one in each of $s_i$ and $[s_i, s^{\varphi}_j]$ are trivial in $\nu(G)$.

\smallskip

\noindent (ii) The $p$-th power of all basic commutators $\{ s_i, [s_i, s^{\varphi}_j] \}$ of weight $< p$ and weight at least one in each of $s_i$ and $[s_i, s^{\varphi}_j]$ are trivial in $\nu(G)$. 

\smallskip

\noindent (iii) $[s_i, s^{\varphi}_j]^p = 1$ in $\nu(G)$.

\smallskip

\noindent (iv) $[P_i(G), G^{\varphi}]$ has exponent at most $p$.

\bigskip

\noindent {\bf Proof.} (i) If $n \leq p$, then $\gamma_{p+1}(\nu(G)) = 1$ (\cite[Corollary 3.2]{roc}) and the statement is immediate. If $n = p+1$, any such basic commutator belong to $\gamma_{2i+1+(p-2)i}(\nu(G)) \subseteq \gamma_{p+2}(\nu(G)) = 1$.

\smallskip

\noindent We prove (ii), (iii) and (iv) jointly using backward induction on $2 \leq i \leq n-1$. 

\smallskip

\noindent First consider $i=n-1$ and let $b$ be a basic commutator in $\{ s_{n-1}, [s_{n-1}, s^{\varphi}_j] \}$ of weight $< p$ and weight at least one in each of $s_{n-1}$ and $[s_{n-1}, s^{\varphi}_j]$. We use induction on weight of $b$ in $\{ s_{n-1}, [s_{n-1}, s^{\varphi}_j] \}$ (which is $\geq 2$) and we denote this by ${\mathrm{wt}}(b)$. 

\smallskip

\noindent If ${\mathrm{wt}}(b) = 2$, then $b$ is either the following element, or its inverse and 
\[
\Big[ [s_{n-1}, s^{\varphi}_j], s_{n-1} \Big] = \Big[ [s_{n-1}, s_j], s^{\varphi}_{n-1} \Big] = 1.
\]
So assume $d := {\mathrm{wt}}(b) \geq 3$ and write $b = [b_1, b_2]$, where $b_1$ and $b_2$ are basic commutators with ${\mathrm{wt}}(b_i) = d_i \geq 1 ~(i=1,2)$ and $d = d_1 + d_2$. Since $b_1 > b_2$ in the lexicographic ordering we have $d_1 \geq 2$. Then $b_1$ is either trivial, or it contains a sub-commutator of the form $[[s_{n-1}, s^{\varphi}_j], s_{n-1}]$ which is a trivial element in $\nu(G)$. Hence $b = 1$ in $\nu(G)$. This proves (ii) for $i={n-1}$. 

\smallskip

\noindent Now we consider (iii) for $i=n-1$. Using Proposition \ref{comm-collection}, we have
\[
1 = [s^p_{n-1}, s^{\varphi}_j] \equiv [s_{n-1}, s^{\varphi}_j]^p \big[ [s_{n-1}, s^{\varphi}_j], s_{n-1} \big]^{p \choose 2} \dotsc \big[ [s_{n-1}, s^{\varphi}_j],_{p-1} s_{n-1} \big] ~{\mathrm{mod}}~ K(s_{n-1}, [s_{n-1}, s^{\varphi}_j]).
\]
Using (i) and (ii) for $i=n-1$, $K(s_{n-1}, [s_{n-1}, s^{\varphi}_j]) = 1$ and using Lemma \ref{exponent-1-phi}, all elements of the right end except the first one are trivial. This proves (iii) for $i=n-1$.

\smallskip

\noindent Next we consider (iv) for $i=n-1$. From Lemma \ref{section-lemma}(i) $[P_{n-1}(G), G^\varphi]$ is an abelian group generated by the elements $[s_{n-1}, s^{\varphi}_j] (0 \leq j \leq n-1)$ each of whose $p$-th power is trivial. This proves (iv) for $i=n-1$.

\smallskip

\noindent Now we assume (ii), (iii) and (iv) for all $k \in \{ i, i+1, \dotsc, n-1 \}$ for some $i \geq 3$. We want to prove (ii), (iii) and (iv) for $k=i-1$.  

\smallskip

\noindent Let $b$ be a basic commutator in $\{ s_{i-1}, [s_{i-1}, s^{\varphi}_j] \}$ of weight $< p$ and weight at least one in each of $s_{i-1}$ and $[s_{i-1}, s^{\varphi}_j]$. If ${\mathrm{wt}}(b) = 2$, then $b$ is either the following element 
\[
\Big[ [s_{i-1}, s^{\varphi}_j], s_{i-1} \Big] = \Big[ [s_{i-1}, s_j], s^{\varphi}_{i-1} \Big],
\]
or its inverse. This element belongs to $[P_i(G), G^{\varphi}]$ and hence its $p$-th power is trivial by induction hypothesis to (iv). 

\smallskip

\noindent If $d := {\mathrm{wt}}(b) \geq 3$ we write $b = [b_1, b_2]$, where $b_1$ and $b_2$ are basic commutators with ${\mathrm{wt}}(b_i) = d_i$. As argued earlier, $d_1 \geq 2$ and hence $b_1$ is either trivial, or it contains a sub-commutator of the form $[[s_{i-1}, s^{\varphi}_j], s_{i-1}]$ which belongs to $[P_i(G), G^{\varphi}]$. Since $[P_i(G), G^{\varphi}]$ is a normal subgroup of $\nu(G)$, it follows that $b \in [P_i(G), G^{\varphi}]$ and consequently $b^p = 1$. This proves (ii) for $i-1$.

\bigskip

\noindent Now we consider (iii) for $i-1$. Using Proposition \ref{comm-collection}, we have
\[
1 = [s^p_{i-1}, s^{\varphi}_j] \equiv [s_{i-1}, s^{\varphi}_j]^p \big[ [s_{i-1}, s^{\varphi}_j], s_{i-1} \big]^{p \choose 2} \dotsc \big[ [s_{i-1}, s^{\varphi}_j],_{p-1} s_{i-1} \big] ~{\mathrm{mod}}~ K(s_{i-1}, [s_{i-1}, s^{\varphi}_j]).
\]
Using (i) and (ii) for $i-1$, we have $K(s_{i-1}, [s_{i-1}, s^{\varphi}_j]) = 1$ and all elements of the right end except the first and the last one are trivial. The last element belong to $\gamma_{(i-1)(p-1)+i}(\nu(G)) \subseteq \gamma_{p+2}(\nu(G))$. Hence $[s_{i-1}, s^{\varphi}_j]^p = 1$. This proves (iii) for $i-1$.

\bigskip

\noindent Now we prove that ${\mathrm{exp}}([P_{i-1}(G), G^{\varphi}]) \leq p$. From Lemma \ref{section-lemma} any element $\xi \in [P_{i-1}(G), G^{\varphi}]$ can be written as $\xi = uw$, where  
\[
u = [s_{i-1}, s^{\varphi}_0]^{\lambda_0} [s_{i-1}, s^{\varphi}_1]^{\lambda_1} \dotsc [s_{i-1}, s^{\varphi}_{n-1}]^{\lambda_{n-1}},
\]
and $w \in [P_i(G), G^{\varphi}]$. From (iii) above applied to $i-1$, we have $0 \leq \lambda_j \leq p-1$ for every $j$. 

\smallskip

\noindent Now, 
\[
(uw)^p \equiv u^p w^p [w,u]^{p \choose 2} \dotsc [w,_{p-1}u] ~{\mathrm{mod}}~ K(u,w).
\]
Any basic commutator in $u,w$ of weight $< p$ and weight at least $2$ in $w$ belong to $[P_i(G), G^{\varphi}] \unlhd \nu(G)$ and hence its $p$-th power is trivial by induction hypothesis. On the other hand any basic commutator in $u,w$ of weight $\geq p$ and weight at least $2$ in $w$ belong to $\gamma_{p+2}(\nu(G))$ and hence it is trivial. This implies that $K(u, w) = 1$. Next for any $1 \leq t \leq p-2$ we have $[w,_{t} u] \in [P_i(G), G^{\varphi}] \unlhd \nu(G)$ and hence their $p$-th power is trivial. Finally $[w,_{p-1} u] \in \gamma_{i+1+(p-1)i}(\nu(G)) \subseteq \gamma_{p+2}(\nu(G)) = 1$. From the above equation we then obtain $(uw)^p = u^p$. Hence it is enough to show that $u^p = 1$. From Lemma \ref{section-lemma} we have $[[P_{i-1}(G), G^{\varphi}], [P_{i-1}(G), G^{\varphi}]] \subseteq [P_i(G), G^{\varphi}]$. Using an induction on the length $\lambda_0 + \dotsc + \lambda_{n-1}$ and similar arguments as above we then have 
\[
u^p = [s_{i-1}, s^{\varphi}_0]^{p \lambda_0} [s_{i-1}, s^{\varphi}_1]^{p \lambda_1} \dotsc [s_{i-1}, s^{\varphi}_{n-1}]^{p \lambda_{n-1}}.
\]
This is trivial from (iii) applied to $i-1$. This proves (iv) for $i-1$. Hence by induction, the lemma follows. \QED

\bigskip

\subsection{Corollary}\label{cor-gen-at-least-wt2} For any $2 \leq i \leq n-1$ we have $[P_i(G), G^{\varphi}] = [P_i(G)^{\varphi}, G]$, and consequently we have ${\mathrm{exp}}([P_i(G)^{\varphi}, G]) \leq p$.

\bigskip

\noindent {\bf Proof.} Since $i \geq 2$, the first statement is immediate from Lemma \ref{exponent-1-phi} and \ref{section-lemma} using an induction on $i$. The consequence follows from Theorem \ref{gen-at-least-wt2}. \QED

\bigskip

\subsection{Proposition}\label{gen-at-least-wt1} Let $p$ be an odd prime and $G$ be a finite $p$-group of maximal class of order $p^n ~(4 \leq n \leq p+1)$ and exponent $p^2$. Then:

\smallskip

\noindent (i) $[s_1, s^{\varphi}_0]^{p^2} = 1 = [s_0, s^{\varphi}_1]^{p^2}$.

\smallskip

\noindent (ii) $([s_0, s^{\varphi}_1][s_1, s^{\varphi}_0])^p = 1$.

\bigskip

\noindent {\bf Proof.} (i) From Lemma \ref{exponent-Pi-subgroups}, and since $\gamma_{p+2}(\nu(G)) = 1$, we have 
\[
1 = [s^{p^2}_1, s^{\varphi}_0] \equiv [s_1, s^{\varphi}_0]^{p^2} \big[ [s_1, s^{\varphi}_0], s_1 \big]^{p^2 \choose 2} \dotsc \big[ [s_1, s^{\varphi}_0],_{p-1} s_1 \big]^{p^2 \choose p} ~{\mathrm{mod}}~ K \big( [s_1, s^{\varphi}_0], s_1 \big).
\]
Now since $[P_2(G), G^{\varphi}] \unlhd \nu(G)$ and ${\mathrm{exp}}([P_2(G), G^{\varphi}]) \leq p$, we have $K \big( [s_1, s^{\varphi}_0], s_1 \big) = 1$ and consequently, $[s_1, s^{\varphi}_0]^{p^2} = 1$. The other equality can be proved similarly. 

\smallskip

\noindent (ii) We set $a := [s_0, s^{\varphi}_1], b := [s_1, s^{\varphi}_0]$. Then we have
\[
(ab)^p \equiv a^p b^p [b,a]^{p \choose 2} \dotsc [b,_{p-1} a] ~{\mathrm{mod}}~ K(a, b).
\]

\noindent Any basic commutator of weight $\geq p$ in $a, b$ with weight at least $2$ in $b$, belong to $\gamma_{4+2(p-2)}(\nu(G)) \subseteq \gamma_{p+2}(\nu(G)) = 1$. Arguing in a similar way, we obtain $[b,_{p-1} a] = 1$. 

\smallskip

\noindent Next, any basic commutator of weight $< p$ and weight at least $2$ in $b$ contains a sub-commutator of the form (or its inverse)
\begin{eqnarray*}
\Big[ [s_1, s^{\varphi}_0], [s_0, s^{\varphi}_1] \Big] & = & \Big[ [s_1, s_0], [s_0, s^{\varphi}_1] \Big] = \Big[ [s_0, s^{\varphi}_1], [s_1, s_0] \Big]^{-1} \\
	& = & \Big[ [s_0, s_1], [s_1, s_0]^{\varphi} \Big]^{-1}, 
\end{eqnarray*}
which belong to $[P_2(G), G^{\varphi}]$. Since $[P_2(G), G^{\varphi}] \unlhd \nu(G)$, it follows that such basic commutator belong to $[P_2(G), G^{\varphi}]$. From Lemma \ref{gen-at-least-wt2}, it follows that the $p$-th power of such basic commutator is trivial. This implies that $K(a,b) = 1$. Arguing in a similar way, we have $[b,_{t} a]^{p \choose t+1} = 1$ for all $1 \leq t \leq p-2$. Thus we get the equation
\[
\Big( [s_0, s^{\varphi}_1][s_1, s^{\varphi}_0] \Big)^p = [s_0, s^{\varphi}_1]^p [s_1, s^{\varphi}_0]^p.
\]

\smallskip

\noindent Since $G/{P_{n-1}(G)}$ has exponent $p$, we may assume that $s^p_1 = s^{\lambda}_{n-1}$ for some $0 \leq \lambda \leq p-1$. Since $s_{n-1}$ is central in $G$, using arguments identical as in the proof of Lemma \ref{section-lemma} (i) we have 
\[
[s_{n-1}, s^{\varphi}_0]^{\lambda} = [s^{\lambda}_{n-1}, s^{\varphi}_0] = [s^p_1, s^{\varphi}_0].
\]
In another way, we have
\[
[s^p_1, s^{\varphi}_0] = [s_1, s^{\varphi}_0]^p \Big[ [s_1, s^{\varphi}_0], s_1 \Big]^{p \choose 2} \dotsc \Big[ [s_1, s^{\varphi}_0],_{p-1} s_1 \Big] ~{\mathrm{mod}}~ K(s_1, [s_1, s^{\varphi}_0]).
\]

\noindent Any basic commutator of weight $\geq p$ in $s_1, [s_1, s^{\varphi}_0]$ with weight at least $2$ in $[s_1, s^{\varphi}_0]$ belong to $\gamma_{4 + (p-2)}(\nu(G)) = 1$. Next, any non-trivial basic commutator of weight $<p$ in $s_1, [s_1, s^{\varphi}_0]$ with weight at least $2$ in $[s_1, s^{\varphi}_0]$ has a sub-commutator of the form (or its inverse)
\[
\Big[ [s_1, s^{\varphi}_0], s_1 \Big] = \Big[ [s_1, s_0], s^{\varphi}_1 \Big],
\] 
which belong to $[P_2(G), G^{\varphi}] \unlhd \nu(G)$, and hence its $p$-th power is trivial. Finally,
\[
\Big[ [s_1, s^{\varphi}_0],_{p-1} s_1 \Big] = \Big[ [[s_1, s_0],_{p-2} s_1], s^{\varphi}_1 \Big].
\]
This implies,
\[
[s_1, s^{\varphi}_0]^p = [s_{n-1}, s^{\varphi}_0]^{\lambda} \Big[ [[s_1, s_0],_{p-2} s_1], s^{\varphi}_1 \Big]^{-1}.
\]
Using similar arguments, we get
\[
[s^{\varphi}_1, s_0]^p = [s^{\varphi}_{n-1}, s_0]^{\lambda} \Big[ [[s_1, s_0],_{p-2} s_1], s^{\varphi}_1 \Big]^{-1}.
\]
From this, we have
\[
[s_1, s^{\varphi}_0]^p [s_0, s^{\varphi}_1]^p = [s_{n-1}, s^{\varphi}_0]^{\lambda} [s^{\varphi}_{n-1}, s_0]^{-\lambda} = 1
\]
using Lemma \ref{exponent-1-phi}. \QED

\bigskip

\noindent We are now ready to prove our final result. 

\bigskip

\subsection{Theorem}\label{main-thm-Schur-mult} Let $G$ be a finite $p$-group of maximal class and order $p^n ~(4 \leq n \leq p+1)$. Then, ${\mathrm{M}}(G)$ is elementary abelian. 

\bigskip

\noindent {\bf Proof.} Using \cite[Lemma 9]{bet}, it follows that $M(G) \neq 1$. If ${\mathrm{exp}}(G) = p$, then the statement follows from \cite[Theorem 1.4]{mor}. So we assume that ${\mathrm{exp}}(G) = p^2$. 

\smallskip

\noindent Let $\alpha \in {\mathrm{M}}(G)$ and we write 
\[
\alpha = \omega_0 \omega_1 \dotsc \omega_{n-1},
\]
where 
\[
\omega_i = \prod_{0 \leq k \leq n-1, k \neq i} [[s_i, s^{\varphi}_k]]^{\lambda_{ik}}. 
\]
The indices satisfy $0 \leq \lambda_{01}, \lambda_{10} \leq p^2 - 1$ and $0 \leq \lambda_{ik} \leq p-1$ if $(i,k) \neq (0,1), (1,0)$. Now notice that $\Psi(\omega_0) \equiv s^{-\lambda_{01}}_2$ mod $P_3(G)$, $\Psi(\omega_1) \equiv s^{\lambda_{10}}_2$ mod $P_3(G)$ and $\Psi(\omega_2 \dotsc \omega_p) \in P_3(G)$. Since $\Psi(\alpha) = 1$, it follows that $\lambda_{10} \equiv \lambda_{01}$ mod $p$. 

\smallskip

\noindent We write
\[
\omega_0 \omega_1 = [[s_0, s^{\varphi}_1]]^{\lambda_{01}} [[s_1, s^{\varphi}_0]]^{\lambda_{10}} \Big( \prod_{k=2}^{n-1} [[s_0, s^{\varphi}_k]]^{\lambda_{0k}} \Big)^{[[s_1, s^{\varphi}_0]]^{\lambda_{10}}} \Big( \prod_{k=2}^{n-1} [[s_1, s^{\varphi}_k]]^{\lambda_{1k}} \Big).
\]

\noindent Since $\lambda_{10} \equiv \lambda_{01}$ mod $p$, we have $[[s_0, s^{\varphi}_1]]^{\lambda_{01}} [[s_1, s^{\varphi}_0]]^{\lambda_{10}} \in {\mathrm{M}}(G)$. This implies that 
\[
\Big( \prod_{k=2}^{n-1} [[s_0, s^{\varphi}_k]]^{\lambda_{0k}} \Big)^{[[s_1, s^{\varphi}_0]]^{\lambda_{10}}} \Big( \prod_{k=2}^{n-1} [[s_1, s^{\varphi}_k]]^{\lambda_{1k}} \Big) \omega_3 \dotsc \omega_p \in {\mathrm{M}}(G).
\]
But this last element is an image of an element of $[P_2(G), G^{\varphi}] = [P_2(G)^{\varphi}, G]$ which has order at most $p$. Hence it is enough to show that $[[s_0, s^{\varphi}_1]]^{\lambda_{01}} [[s_1, s^{\varphi}_0]]^{\lambda_{10}}$ has order $p$ in ${\mathrm{M}}(G)$. 

\smallskip

\noindent Now we denote $a := [[s_0, s^{\varphi}_1]], b := [[s_1, s^{\varphi}_0]]$ and write $\lambda := \lambda_{01}, \lambda_{10} - \lambda_{01} = p\mu$ for some $\mu \in {\mathbb{Z}}$. Then,
\[
[[s_0, s^{\varphi}_1]]^{\lambda_{01}} [[s_1, s^{\varphi}_0]]^{\lambda_{10}} = (ab) (ab)^b \dotsc (ab)^{b^{\lambda - 1}} (b^p)^{\mu}.
\] 
Now $ab, b^p \in M(G)$ both of which has order $p$ from Proposition \ref{gen-at-least-wt1}. Since ${\mathrm{M}}(G)$ is an abelian normal subgroup of $[[G, G^{\varphi}]]$, the statement follows. \QED

\bigskip

\noindent {\bf Acknowledgement.} The authors would like to thank Professor Bettina Eick and Professor Michael VaughanLee for many helpful comments and encouragement. The research of the author Renu Joshi is supported by PMRF (India) fellowship. 

\bigskip
%\newpage

%%%%%%%%%%%%%%%%%%%%%%%%%%%%%%%%%%%%%%%%%%%%%%%%%%%%%%%%%%%%%%%%%%%%%%%%% 
 
\bibliographystyle{plain} 
\bibliography{SchurMultMaxClassV2}

\end{document}